\documentclass[preprint,12pt]{elsarticle}

\usepackage{amsmath, amssymb}
\usepackage{hyperref}
\usepackage{float}

\journal{}

\begin{document}

\begin{frontmatter}
\title{Analytical Standard Errors for Exploratory Factor Solutions}
  
\author[affil1]{Xingwei Hu\footnote{Corresponding author}}
\author[affil2]{Caihong Hu}
\author[affil3]{Cheng-Kuang Wu}
\affiliation[affil1]{
 organization={International Monetary Fund}, 
 city={Washington},   state={DC}, postcode={20431}, country={USA; xhu@imf.org}
}
\affiliation[affil2]{
organization={Zhangjiajie College}, 
city={Zhangjiajie}, state={Hunan}, postcode={427000}, country={China; 2006025@zjjc.edu.cn}
}
\affiliation[affil3]{
organization={Zhejiang Industry and Trade Vocational College}, 
city={Wenzhou}, state={Zhejiang}, postcode={325700}, country={China; shapleyvalue@hotmail.com}
}

\begin{abstract}
Inference for factor models is often hampered by the lack of tractable and accurate variance estimates, which can materially distort downstream analyses.  
In practice, uncertainty in the residual covariance matrix is frequently either ignored or addressed through computationally intensive resampling methods that tend to be unstable. 
This paper develops a unified analytical framework for inference in exploratory factor analysis under several widely used extraction rules, including least-squares, principal-factor, iterative principal-component, alpha, and image factoring. 
By treating these estimators as implicitly defined functions of the sample covariance matrix, we derive closed-form Jacobians that translate perturbations in the covariance matrix into changes in the resulting factor solutions. 
Combined with the delta method and consistent estimators of the sample covariance matrix, the proposed approach yields standard errors that are straightforward to compute and remain valid under non-Gaussianity, heteroskedasticity, and serial or cross-sectional dependence. 
Simulation evidence confirms that the analytical standard errors accurately capture finite-sample variability while avoiding both the instability of bootstrap procedures and the restrictive assumptions underlying Fisher information-based inference. 
An application to a factor-augmented structural vector autoregressive (SVAR) model further demonstrates how accounting for this source of uncertainty can substantially affect impulse-response inference. 
Taken together, the results provide a practical and general tool for propagating estimation uncertainty in settings where factor extraction serves as an intermediate step.
\end{abstract}


\end{frontmatter}

\section{Introduction}\label{sect:intro}

Factor-analytic representations provide a flexible framework for dimension reduction and structural modeling. 
A useful way to build intuition is to start from the residual covariance matrix.
In multivariate models, residuals are typically not independent; rather, they tend to co-move, reflecting underlying common influences that are not directly observed. 
Factor analysis formalizes this idea by representing the covariance matrix as the combined effect of a small number of shared latent forces and many variable-specific disturbances. 
Within this framework, co-movement across variables is attributed to a limited set of latent factors, while the remaining variation is treated as idiosyncratic. 

Formally, for a $p\times p$ variance or correlation matrix $\Sigma$ of residuals, factor analysis imposes the decomposition
\begin{equation}\label{eq:model}
\Sigma = \Lambda \Lambda^\top + \Psi,
\end{equation}
where $\Lambda \in \mathbb{R}^{p\times k}$ denotes the loading matrix with $k<p$, and $\Psi$ is a diagonal matrix of idiosyncratic variances, commonly referred to as uniquenesses. 
This structure yields a parsimonious reduced-form representation of dependence driven by a small number of common or clustered shocks.
In econometrics, psychometrics, and statistics, such a representation is especially useful when the observed co-movement among variables is generated by a few latent source of variation.
For example, in a factor-augmented SVAR model, residuals can be interpreted as the sum of common shocks and idiosyncratic disturbances, with factor analysis providing a tractable and economically meaningful way to separate these components.

The sample residual covariance matrix $\widehat \Sigma$ can be estimated in various ways, depending on the underlying multivariate model, the available data, and the parameter estimation method.
As a by-product, model estimation typically yields an estimate $\widehat \Sigma$, but rarely provides its sampling covariance, denoted by $\operatorname{cov}(\widehat \Sigma)$.
When the residuals are multivariate normal, $\widehat \Sigma$ follows a scaled Wishart distribution.
In many applications, however, the observed variables are binary, categorical, ordinal, bounded, or integer-valued.
Moreover, in time-series and panel-data settings, residuals often exhibit serial dependence and heteroskedasticity, making distribution-free or robust covariance-based inference particularly attractive (e.g., \cite{Lin&Chen&Wang2023,Raykov&Penev&Zhang2026}). 

Given an estimate $\widehat \Sigma$, extracting $\widehat \Lambda$ and $\widehat \Psi$ from Eq.~(\ref{eq:model}) is not  merely a descriptive exercise; it is central to structural inference, including impulse response analysis, forecast-error variance decompositions, and related objects.
In addition to likelihood-based approaches, a broad class of non-maximum-likelihood (non-ML) factor-extraction methods---such as principal-factor, least-squares, alpha, and image factor analysis---remain widely used.
These methods are computationally efficient, perform robustly in practice, and rely on relatively weak distributional assumptions.
More generally, when parametric distributional assumptions are overly restrictive, distribution-free covariance estimators can be constructed for a wide range of statistics, including $\widehat \Sigma$.

Although point estimation of factor loadings and uniquenesses is relatively straightforward across extraction methods, rigorous statistical inference remains considerably more challenging. 
In practice, $\operatorname{cov}(\widehat \Sigma)$ is often treated as a nuisance parameter and omitted from standard software output.
Nevertheless, explicitly accounting for this covariance can substantially improve inference by more accurately capturing asymmetry, sampling variability, and finite-sample bias, albeit at some additional computational cost.
This issue is especially important for downstream analyses---such as forecasting and policy analysis---that depend on uncertainty in nuisance parameters.
Ignoring this source of uncertainty can lead to biased standard errors, overly optimistic forecasts, and, ultimately, misleading policy conclusions and decision-making.

In the existing literature, closed-form standard-error formulas for factor loadings and uniquenesses are available almost exclusively for maximum-likelihood (ML) factorization, where results based on the Fisher-information matrix follow directly from likelihood theory. 
However, ML procedures typically rely on strong distributional assumptions that may fail in the presence of heteroskedasticity or autocorrelation.
This motivates the use of heteroskedasticity- and autocorrelation-consistent (HAC) or heteroskedasticity-consistent (HC) estimators for $\Sigma$.  
Moreover, formal empirical distribution function (EDF) tests for the assumed distribution of $\widehat \Sigma$ are typically unavailable, potentially undermining the validity of ML-based extraction.

For non-ML methods, despite their widespread use in applied work and their implementation in standard statistical software, explicit covariance expressions remain largely unavailable.
These methods are often preferred in settings characterized by non-normality, heavy-tailed distributions, or heterogeneous data, where likelihood-based assumptions are less credible. 
When inference is required, it is therefore typically obtained only indirectly, through resampling techniques.
Resampling methods, however, face several challenges beyond their high computational cost.
They require either an explicit characterization of the sampling distribution of $\widehat\Sigma$ or extensive bootstrap resampling of the original data.
Moreover, bootstrapping can introduce an alignment problem when estimating $\Lambda$, because the ordering---and in some cases, the signs---of estimated factors may vary across resamples. 
Such label switching can inflate estimates of $\operatorname{cov}(\widehat \Lambda)$ and, in turn, produce spuriously large standard errors.

This paper addresses these limitations by deriving explicit closed-form expressions for the covariance matrices of the factor-loading and uniqueness estimators associated with several non-ML extraction methods. 
Specifically, we consider least-squares, principal-factor, iterative principal-component, alpha, and image factorization methods. 
Our approach expresses the covariance matrices of these estimators directly in terms of $\operatorname{cov}(\widehat \Sigma)$. 
As a result, the proposed formulas apply under both normal and non-normal sampling schemes and integrate naturally with existing covariance estimators for $\widehat \Sigma$, including HAC and HC estimators.

Methodologically, our analysis relies on multivariate implicit differentiation of the estimating equations that define each factor solution. 
By treating the estimated loading matrix as an implicit differentiable function of the sample covariance matrix, we apply the delta method to derive variance expressions for the corresponding factor estimators. 
This unified framework yields closed-form Jacobians that translate perturbations in $\widehat \Sigma$ into changes in the estimated factor loadings and uniquenesses, thereby producing tractable covariance formulas.
These Jacobians also provide the foundation for delta-method inference in subsequent analyses, including rotated factor solutions
(\cite{Archer&Jennrich1973}, \cite{Jennrich1973}).

To illustrate the implementation, we present two empirical applications: a least-squares factor analysis and a factor-augmented SVAR model.
The SVAR application includes three endogenous time series but only two structural shocks
and explicitly accounts for $\mathrm{cov}(\widehat \Sigma)$, showing that this additional source of uncertainty materially affects the resulting impulse-response analysis.
In addition, we conduct a Monte Carlo study to assess the finite-sample performance of the proposed standard errors for the uniqueness estimators. 
The simulation results show that the empirical standard errors converge to their theoretical counterparts, providing strong support for the analytical derivations.

From a practical perspective, the proposed results offer several advantages. 
First, they provide a transparent method for computing standard errors without relying on simulation-based resampling or likelihood-specific assumptions. 
Second, they are especially well suited to settings in which factor analysis serves as an intermediate step, such as SVARs with latent structural shocks, because they allow uncertainty in factor estimates to be propagated coherently into downstream quantities of interest.  
Without such propagation, uncertainty may be understated, leading to overconfident forecasts and policy conclusions.
Third, because the formulas depend only on $\widehat \Sigma$ and $\operatorname{cov}(\widehat \Sigma)$, they accommodate a wide range of data-generating processes, including non-Gaussian and high-frequency environments.

The remainder of the paper is organized as follows. 
Section~\ref{sect:review} reviews exploratory factor analysis and the implicit-differentiation approach used throughout the paper.
Section~\ref{sect:IPC_PFA}---\ref{sect:ImageFA} derive the covariance matrices of the factor loadings for iterative principal-component, principal-factor, least-squares, alpha, and image factor analysis. 
Section~\ref{sect:uniqueness} presents closed-form expressions for the standard errors of the uniqueness estimators. 
Section~\ref{sect:Example} reports the two empirical applications and the simulation study.
Section~\ref{sect:diss} concludes and outlines several avenues for future research.

\section{Framework and Methodology}\label{sect:review}

In macroeconomics and finance, the covariance structure in Eq.~(\ref{eq:model}) underpins dynamic factor models, structural vector autoregressions (SVARs), and variance-covariance decompositions of large systems, in which a small number of latent shocks drives co-movement across a high-dimensional set of observables. 
A central objective is to decompose the residuals into two distinct sources of unobserved variation---common and idiosyncratic shocks---each carrying its own economic interpretation and analytical role.

\subsection{Structural Modeling by Exploratory Factors}

For a $p$-dimensional random vector $Y = (X_1, \dots, X_p)^\top$,
classical exploratory factor analysis posits the existence of $k$ mutually independent common factors that underlie the observed variables in $Y$. 
Collecting these latent factors in the random vector $F$, the relationship between the observed variables and the latent factors is specified as
\begin{equation}\label{eq:CFA}
Y = c + \Lambda F + \epsilon,
\end{equation}
where $c$ is a constant vector and $\epsilon$ is a $p$-dimensional vector of mutually independent idiosyncratic disturbances.
Under the standard normalization that the common factors have unit variance, Eq.~(\ref{eq:CFA}) implies
the covariance structure of $Y$ given in Eq.~(\ref{eq:model}).
In particular, the diagonal elements of $\Psi$ correspond to the variances of the idiosyncratic disturbances, commonly referred to as the \textit{uniquenesses}.

Eq.~(\ref{eq:model}) further implies that the off-diagonal elements of the covariance matrix of $Y$ are determined entirely by the factor-loading matrix; the common factors thus drive the co-movement of $Y$.  
By contrast, the diagonal elements reflect contributions from both the common and idiosyncratic components: the common factors capture variation shared across the variables in $Y$, whereas the idiosyncratic disturbances account for variable-specific noise.
Let $\lambda_{it}$ denote the $(i,t)$-th entry of $\Lambda$.
The $i$-th diagonal element of $\Lambda \Lambda^\top$ is then
\begin{equation}\label{eq:h_i}
h_i = \sum\limits_{t=1}^k \lambda_{it}^2,
\end{equation}
which measures the portion of the variance of $X_i$ explained by the common factors and is commonly referred to as the \textit{communality} of $X_i$.

The matrix $\Sigma$ need not correspond to the variance-covariance matrix of a single set of observed variables.
More generally, it may characterize the covariance structure of reduced-form residuals arising from a system of equations. 
As an illustration, consider the following factor-augmented structural vector autoregressive (SVAR) model:
\begin{equation}\label{eq:SVAR}
Y_t = c + \sum_{s=1}^{\zeta} A_s Y_{t-s} + \Lambda F_t + \epsilon_t.
\end{equation}
Here, $Y_t$ is a $p$-dimensional vector of observed time series at time $t$, $\zeta$ denotes the lag length, $F_t$ collects $k$ common structural shocks,
and $\epsilon_t$ contains $p$ idiosyncratic shocks. 
The composite residual $u_t \equiv \Lambda F_t +\epsilon_t$ therefore exhibits a two-component error structure, and $\Sigma$ can be interpreted as the covariance matrix of $u_t$.

For panel data, Eq.~(\ref{eq:SVAR}) extends naturally to
$$
Y_{it} = c_i + \sum_{s=1}^{\zeta} A_s Y_{i,t-s} + \Lambda_i F_t + \epsilon_{it},
$$
where $\Lambda_i$ is a cross-section-specific loading matrix and $c_i$ is a cross-section-specific intercept,
while additional fixed or random effects are omitted for simplicity.
The common factors $F_t$ capture both cross-sectional and cross-variable interdependence in $Y$.
Estimation of $A_s$ and $c_i$ generally requires some structural assumption on the composite residuals.
For least-squares estimation, for instance, $\mathrm{cov}\left(\Lambda_i F_t + \epsilon_{it}\right)$ is assumed to be constant across $i$ and $t$, while allowing for heteroskedasticity in $\epsilon_{it}$.

Much of the SVAR literature ignores idiosyncratic shocks, effectively imposing the restriction $k=p$.
Explicitly specifying idiosyncratic shocks, however, reduces both the number of structural shocks and the number of identifying restrictions required. 
From a computational standpoint, embedding an error-component structure---such as a common-factor specification---within an SVAR is straightforward: the parameters $c$ and $A_s$ in Eq.~(\ref{eq:SVAR}) can be estimated by least squares under the assumption that the variance-covariance matrix of $\Lambda F_t + \epsilon_t$ is time-invariant.
Despite the availability of such results, uncertainty in $\widehat \Sigma$ is rarely used when defining structural shocks and is not fully propagated into subsequent inference, such as impulse-response analysis.
For example, Cholesky shocks are typically obtained by applying a Cholesky decomposition to $\widehat \Sigma$ without accounting for $\operatorname{cov}(\widehat \Sigma)$.

This form of dimensionality reduction is particularly valuable in settings characterized by highly interconnected data.
In macroeconomic modeling, structural shocks are typically assumed to be mutually independent. 
Yet, although an economy can be examined from many angles---yielding a large set of interrelated observable time series---macroeconomic fluctuations are generally driven by only a small number of underlying shocks, such as demand, 
supply, geopolitical, and technological shocks.
It is also unlikely that the unexpected residual in any given observable can be attributed solely to this limited set of structural shocks, without any influence from other relevant disturbances.
Consequently, imposing $p=k$ may force researchers to exclude informative data series, giving rise to omitted-variable bias and potentially introducing heteroskedasticity and serial autocorrelation into the residuals.

\subsection{Multivariate Implicit Differentiation}
A wide range of factor-extraction methods has been proposed in the literature (e.g., \cite{Harman1976,Lawley&Maxwell1971,SAS2015}), including principal-factor, maximum-likelihood, and least-squares approaches.
These methods aim to recover the factor model in Eq.~(\ref{eq:model}) from an estimated sample covariance (or correlation) matrix $\widehat \Sigma$.
From a practical standpoint, however, it is important not only to obtain point estimates of the factor loadings $\widehat \Lambda$ and uniquenesses $\widehat \Psi$, but also to quantify their precision and to translate them into substantively meaningful applications---such as factor rotation.
Achieving these objectives requires both the covariance matrices of the estimators and 
a clear understanding of how those covariances propagate through the chosen rotation and related criteria.

The covariance of the loading estimator $\widehat \Lambda$ depends jointly on the factor-extraction method and on
the covariance of $\widehat \Sigma$. 
Let the population loading matrix $\Lambda$ be defined as a function of the population covariance matrix $\Sigma$, namely, $$
\Lambda = f(\Sigma),
$$
with sample counterpart $\widehat \Lambda = f(\widehat \Sigma)$ under a given extraction rule $f$.
Assuming that $f$ is differentiable at $\Sigma$, by the delta method, small permutations in $\Sigma$ propagate to $\Lambda$ through the Jacobian of $f$, yielding
$$
d\, \Lambda = \frac{\partial f}{\partial \Sigma} d\, \Sigma.
$$
In this setting, $d \Lambda \approx \widehat \Lambda - \Lambda$ and $d\Sigma \approx \widehat \Sigma - \Sigma$ represent deviations from the unobserved population quantities $\Lambda$ and $\Sigma$, respectively. 
In particular, when $\widehat \Lambda$ and $\widehat \Sigma$ are unbiased estimators of $\Lambda$ and $\Sigma$, the linearization implies
\begin{equation}\label{eq:cov_Lambda_hat}
\operatorname{cov}\left(\widehat \Lambda\right) 
= 
\left[\frac{\partial f}{\partial \Sigma}\right] \operatorname{cov}\left(\widehat \Sigma\right) 
\left[\frac{\partial f}{\partial \Sigma}\right]^\top.
\end{equation}
When $\widehat \Lambda$ and $\widehat \Sigma$ are only consistent, this equality is interpreted as an asymptotic approximation that becomes exact as the sample size tends to infinity.

Writing $\Sigma = [\sigma_{ij}]$, and letting $\widehat \sigma_{ij}$ denote the sample covariance 
(or correlation) between $X_i$ and $X_j$, Eq.~(\ref{eq:cov_Lambda_hat}) implies that the covariance between entries of $\widehat \Lambda$ can be expressed as 
\begin{equation}\label{eq:cov_Lambda}
\operatorname{cov}\left(\widehat \lambda_{ir}, \widehat\lambda_{js}\right) 
= 
\sum_{m,n,x,y=1}^p \frac{\partial \lambda_{ir}}{\partial \sigma_{mx}}
\operatorname{cov}\left(\widehat \sigma_{mx},\widehat \sigma_{ny}\right) 
\frac{\partial \lambda_{js}}{\partial \sigma_{ny}}.
\end{equation}
Analogously, let $\psi = (\psi_1,\ldots,\psi_p)^\top$ denote the vector of diagonal elements of the uniqueness matrix $\Psi$. 
The covariance of the uniqueness estimators is then
\begin{equation}\label{eq:cov_Psi}
\operatorname{cov}\left(\widehat \psi_i, \widehat \psi_j\right)
= 	
\sum_{m,n,x,y=1}^p \frac{\partial \psi_i}{\partial \sigma_{mx}}
\operatorname{cov}\left(\widehat \sigma_{mx},\widehat \sigma_{ny}\right) 
\frac{\partial \psi_j}{\partial \sigma_{ny}}.
\end{equation}
Similarly, the delta method implies that the cross-covariance between $\widehat \Lambda$ and $\widehat \psi$ takes the form
\begin{equation}\label{eq:cov_Lambda_Psi}
\operatorname{cov}\left(\widehat \Lambda, \widehat \psi \right)
=
\left[\frac{\partial f}{\partial \Sigma}\right] \operatorname{cov}\left(\widehat \Sigma\right) 
\left[\frac{\partial \psi}{\partial \Sigma}\right]^\top.
\end{equation}

In general, however, $\Lambda = f(\Sigma)$ is only implicitly defined by the factor-extraction rule.
Fortunately, an explicit form of $f$ is not required: only the derivatives $\frac{\partial \Lambda}{\partial \Sigma}$
and $\frac{\partial \Psi}{\partial \Sigma}$ enter Eqs.(\ref{eq:cov_Lambda_hat})---(\ref{eq:cov_Lambda_Psi}).
These derivatives can be obtained through implicit differentiation of the identifying condition linking $\Lambda$ and $\Sigma$. 
For instance, suppose the factor-extraction condition takes the form
$$
\Phi (\Lambda, \Sigma) = 0, \quad \mathrm{or}\ \mathrm{equivalently} \quad \Phi (f(\Sigma), \Sigma)=0.
$$
Differentiating either expression with respect to $\Sigma$ yields
$$
\frac{\partial \Phi}{\partial \Lambda} \frac{\partial f}{\partial \Sigma} 
+
\frac{\partial \Phi}{\partial \Sigma}  = 0,
$$
and therefore
$$
\frac{\partial f}{\partial \Sigma} 
=
- \left( \frac{\partial \Phi}{\partial \Lambda} \right)^{-1} \frac{\partial \Phi}{\partial \Sigma}.
$$

In the next five sections, we derive explicit expressions for the differentials $d\Lambda$ and $d \psi$ as linear functions of $d\Sigma$ for several widely used factor-extraction conditions of the form $\Phi=0$. 
Combined with the covariance matrix of $\widehat \Sigma$ (e.g.,\cite{deLeeuw1983, Girshick1939,Hsu1949}), these results yield standard errors for the estimated loadings $\widehat \Lambda$ and uniquenesses $\widehat \psi$, and provide the essential ingredients for delta-method inference on downstream analyses via the chain rule. 
While earlier work has developed asymptotic theory for principal components and standard-error formulas for maximum-likelihood factor analysis (e.g.,\cite{Anderson1963,Girshick1939,Jennrich&Thayer1973,Lawley1967}), our focus is on widely used non-ML procedures as implemented in standard statistical and econometric software. 
Because these methods are pervasive in applications involving discrete or qualitative measurements, high dimensionality, and heavy-tailed data, the absence of closed-form covariance expressions constitutes an important gap that hinders both formal inference and rigorous model comparison.

\section{Iterative Principal Component and Principal Factor Analysis}\label{sect:IPC_PFA}
We begin by introducing notation. 
For the $p\times k$ matrix $\Lambda$, let $\overrightarrow{\Lambda}$ denote
its vectorization obtained by stacking the columns of $\Lambda$ from left to right, that is, 
$$
\overrightarrow{\Lambda} = (\lambda_{11},\ldots,\lambda_{p1},\lambda_{12},\ldots,\lambda_{p2},\ldots,\lambda_{1k},\ldots,\lambda_{pk})^\top.
$$
When no ambiguity arises, we use the same double-subscript notation to refer to elements of $\overrightarrow{\Lambda}$. 
In particular, the $(i,r)$-th element of $\Lambda$, denoted by $\lambda_{ir}$,
corresponds to position $((r-1)p+i)$ in $\overrightarrow{\Lambda}$.
For a square matrix $M$, $\operatorname{diag} (M)$ denotes the vector of diagonal entries of $M$,
and $\operatorname{Diag} (M)$ denotes the diagonal matrix whose diagonal equals $\operatorname{diag} (M)$.
Finally, $\delta_{ij}$ denotes the Kronecker delta, which equals $1$ if $i=j$ and $0$ otherwise.

Both iterative principal-component factor analysis and principal factor analysis construct the loading matrix $\Lambda$ from the leading eigenvectors of $\Sigma - \Psi$:
\begin{equation}\label{eq:SD}
\Sigma - \Psi = \Lambda \Lambda^\top = \sum_{r=1}^k \lambda_r \lambda_r^\top,
\end{equation} 
where $\lambda_r$ denotes the $r$-th column of $\Lambda = [\lambda_1,\lambda_2,\ldots,\lambda_k]$.
These columns are ordered by the eigenvalues of $\Sigma - \Psi$, with $\lambda_r$ associated with the $r$-th largest eigenvalue.
Without loss of generality, assume these eigenvalues are distinct, as is typically the case for a sample variance-covariance matrix.
Under this condition, the columns of $\Lambda$ are mutually orthogonal, so the inner product between any two distinct columns is zero.
Taking diagonal elements of Eq.~(\ref{eq:SD}) then yields 
\begin{equation}\label{eq:L2}
	\sigma_{ii}-\psi_i = \sum_{t=1}^k \lambda_{it}^2, \quad \forall \ i=1,...,p.
\end{equation}

For the sample variance-covariance matrix $\widehat \Sigma$, Eq.~(\ref{eq:SD}) may, in general, contain additional terms of the form $\widehat\lambda_r \widehat \lambda_r^\top$ for $r>k$.
These items are omitted in accordance with the model specification in Eq.~(\ref{eq:model}) and accounted for by sampling errors.

Using Eq.~(\ref{eq:SD}), we obtain
\begin{equation}\label{eq:pfa-se-eigen}
\left(\Sigma-\Psi\right) \lambda_r 
= \Lambda \Lambda^\top \lambda_r 
= \sum\limits_{s=1}^k \lambda_s \lambda_s^\top \lambda_r
= \sum\limits_{s=1}^k \lambda_s \left(\lambda_s^\top \lambda_r\right)
= \left(\lambda_r^\top \lambda_r\right) \lambda_r.
\end{equation}
This identity holds for any $1\le r \le k$. 
Fix indices $1\le i\le p$ and $1\le r \le k$. 
Taking the $i$-th component of both sides of Eq.~(\ref{eq:pfa-se-eigen}) and applying Eq.~(\ref{eq:L2}) yields
$$
\sum_{j\not =i} \sigma_{ij}\lambda_{jr} + \lambda_{ir}\sum_{t=1}^k \lambda_{it}^2 
= \lambda_{ir}\sum_{j=1}^p \lambda_{jr}^2.
$$
Taking differentials on both sides gives
$$
\resizebox{.98\textwidth}{!}{$
\sum\limits_{j\not =i} \left(\sigma_{ij} d\lambda_{jr} + \lambda_{jr} d \sigma_{ij}\right) 
+ \left(\sum\limits_{t=1}^k \lambda_{it}^2\right) d \lambda_{ir}+ 2 \lambda_{ir}\sum\limits_{t=1}^k \lambda_{it} d \lambda_{it}
= 
\left(\sum\limits_{j=1}^p \lambda_{jr}^2\right) d \lambda_{ir} + 2 \lambda_{ir}\sum\limits_{j=1}^p \lambda_{jr} d \lambda_{jr}.
$}
$$
After algebraic simplification, this expression reduces to
$$
\resizebox{.98\textwidth}{!}{$
\sum\limits_{j\not =i} \left( \sigma_{ij}-2\lambda_{ir}\lambda_{jr} \right) d\lambda_{jr}
+\left( \sum\limits_{t=1}^k \lambda_{it}^2-\sum\limits_{j=1}^p \lambda_{jr}^2 \right) d\lambda_{ir}
+ 2\lambda_{ir} \sum\limits_{t\not =r}\lambda_{it}d\lambda_{it}
= 
- \sum\limits_{j\not = i} \lambda_{jr} d\sigma_{ij}.
$}
$$

Collecting the linear equations derived above, we can write the system compactly in matrix form as
$$
A d\overrightarrow{\Lambda} = B d\overrightarrow{\Sigma},
$$ 
where the coefficient of $d\lambda_{jt}$ in the $(i,r)$-th equation is 
$$
A_{ir,jt} = \left \{
\begin{array}{ll}
	\sum\limits_{s=1}^k \lambda_{is}^2-\sum\limits_{s=1}^p \lambda_{sr}^2,	\ & \ \mathrm{if}\ j=i \ \mathrm{and}\ t=r,\\
	2\lambda_{ir}\lambda_{it},	            \ & \ \mathrm{if}\ j=i \ \mathrm{but}\ t\not =r,\\
	\sigma_{ij}-2\lambda_{ir}\lambda_{jr},    \ & \ \mathrm{if}\ j\not =i \ \mathrm{but}\ t =r,\\
	0,				                  \ & \ \mathrm{otherwise}.\\
\end{array}\right .
$$
This holds for all $1\le i,j\le p$ and $1\le r,t\le k$.
Similarly, the coefficient of $d\sigma_{xy}$ in the same equation is
$$
B_{ir,xy} = -\delta_{ix}(1-\delta_{xy})\lambda_{yr},
\quad
\forall\ 1\le i, x,y\le p, \ 1\le r\le k. 
$$

Therefore, 
$$
d \overrightarrow{\Lambda} = A^{-1} B d\overrightarrow{\Sigma}.
$$ 
By the delta method, the covariance matrix of $\widehat{\Lambda}$ can be estimated as
$$
\operatorname{cov}\left( \overrightarrow{\widehat{\Lambda}} \right) 
=
\left(A^{-1} B\right) \operatorname{cov}\left( \overrightarrow{\widehat{\Sigma}} \right) \left(A^{-1} B\right)^\top
$$ 
where $A$ and $B$ are evaluated at $\widehat \Lambda$ (rather than $\Lambda$).

\section{Least-Squares Factor Analysis}\label{sect:LSFA}
In least-squares factor analysis (e.g., \cite{Harman1976}), the loading matrix $\Lambda$ is chosen to minimize the squared off-diagonal discrepancies in the residual covariance matrix $\Sigma- \Lambda  \Lambda^\top$. 
Specifically, define the objective function
\begin{equation}\label{eq:LS_objective}
	\mathrm{g}(\Lambda) =  \sum_{1\le y < x \le p} \left(\sigma_{xy} - \sum_{t=1}^k \lambda_{xt} \lambda_{yt} \right)^2,
\end{equation}
that is, the sum of squares of the elements in the strictly lower-triangular part of $\Sigma-\Lambda\Lambda^\top$.
When the number of free parameters in $\Lambda$ is sufficient to drive $\mathrm{g}(\Lambda)$ to zero (i.e., to match all off-diagonal covariances exactly), the criterion can be extended to include the diagonal elements of $\Sigma$ while maintaining $\mathrm{g}(\Lambda)=0$. 
Doing so further increases the portion of total variance explained by the common factors. 

The first-order optimality conditions impose $p\times k$ restrictions on $\Lambda$.
Specifically, we require
$$
\frac{\partial \mathrm{g} (\Lambda)}{\partial \lambda_{ir}} = 0
$$ 
for all $1\le i \le p$ and $1\le r \le k$.
A straightforward differentiation yields
$$
\begin{array}{rcl}
	-\frac{1}{2} \frac{\partial \mathrm{g}}{\partial \lambda_{ir}}(\Lambda)
&=&
\sum\limits_{_{1\le y<x\le p}} \left(\sigma_{xy} 
	- \sum\limits_{t=1}^k \lambda_{xt} \lambda_{yt} \right)
	\sum\limits_{t=1}^k (\delta_{ix} \delta_{rt} \lambda_{yt} + \delta_{iy} \delta_{rt} \lambda_{xt}) \\
&=&
\sum\limits_{_{1\le y <x \le p}} (\delta_{ix} \lambda_{yr} + \delta_{iy} \lambda_{xr})
	\left(\sigma_{xy} - \sum\limits_{t=1}^k \lambda_{xt} \lambda_{yt} \right) \\
&=&
\sum\limits_{_{1\le y<i}} \lambda_{yr} \left(\sigma_{iy} - \sum\limits_{t=1}^k \lambda_{it} \lambda_{yt} \right) 
	+ \sum\limits_{_{i<x\le p}} \lambda_{xr} \left(\sigma_{xi} - \sum\limits_{t=1}^k \lambda_{xt} \lambda_{it} \right) \\
&=&
\sum\limits_{_{j\not =i}} \lambda_{jr} \left(\sigma_{ij} 
	- \sum\limits_{t=1}^k \lambda_{it} \lambda_{jt} \right).
\end{array}
$$
Therefore, for any fixed $1\le i\le p$ and $1\le r\le k$, the first-order condition reduces to
\begin{equation}\label{eq:lse-equal}
\sum_{j\not = i} \lambda_{jr}\sigma_{ij} = \sum_{j\not = i} \lambda_{jr} \sum_{t=1}^k \lambda_{it} \lambda_{jt}.
\end{equation}

The differential form of Eq.~(\ref{eq:lse-equal}) is 
$$
\resizebox{.98\textwidth}{!}{$
\sum\limits_{j\not = i} \left(\lambda_{jr} d\sigma_{ij}+ \sigma_{ij} d \lambda_{jr} \right) 
= 
\sum\limits_{j\not = i} \left( \sum\limits_{t=1}^k \lambda_{it} \lambda_{jt}\right) d\lambda_{jr}
+
\sum\limits_{j\not = i} \lambda_{jr} \sum\limits_{t=1}^k \left( \lambda_{jt} d \lambda_{it} + \lambda_{it} d \lambda_{jt}\right).
$}
$$
Collecting like terms and simplifying yield
$$
\resizebox{.98\textwidth}{!}{$
\sum\limits_{j\not = i} \lambda_{jr} d\sigma_{ij}
= 
\sum\limits_{j\not = i}\left(\sum\limits_{t=1}^k \lambda_{it} \lambda_{jt}-\sigma_{ij}+\lambda_{jr}\lambda_{ir}\right) d\lambda_{jr} 
+ \sum\limits_{t=1}^k \left(\sum\limits_{j\not = i} \lambda_{jr} \lambda_{jt}\right) d\lambda_{it}
+ \sum\limits_{j\not =i}\sum\limits_{t\not =r} \lambda_{jr} \lambda_{it}d \lambda_{jt},
$}
$$
which completes the differential characterization of the first-order condition.

We next rewrite the resulting system in matrix form as
$$
-B d\overrightarrow{\Sigma}  = D d \overrightarrow{\Lambda},
$$
where the coefficient of $d\lambda_{jt}$ in the $(i,r)$-th equation is defined by
$$
D_{ir,jt} = \left \{
\begin{array}{ll}
\sum\limits_{z\not = i}\lambda_{zr}\lambda_{zt},& \quad \mathrm{if}\ j=i,\\
\sum\limits_{z=1}^k \lambda_{iz}\lambda_{jz}-\sigma_{ij}+\lambda_{ir}\lambda_{jr},
& \quad \mathrm{if}\ j\not =i\ \mathrm{and} \ t=r,\\
\lambda_{jr}\lambda_{it},	   & \quad \mathrm{if}\ j\not =i\ \mathrm{and} \ t\not =r.\\
\end{array} \right .
$$ 
These definitions apply for all $1\le i,j\le p$ and $1\le r, t \le k$. 
Solving for $d\overrightarrow{\Lambda}$ gives
$$
d\overrightarrow{\Lambda} = - D^{-1} B d\overrightarrow{\Sigma}.
$$

\section{Alpha Factor Analysis}\label{sect:ALphaFA}
Let $H$ be the diagonal matrix with $i$-th diagonal element $\sqrt{h_i}$. 
Clearly, $h_i>0$; otherwise, Eq.~(\ref{eq:h_i}) would imply that the $i$-th row of the loading matrix $\Lambda$ is the zero vector, making the associated factor redundant. 
Moreover, $H^{-1} (\Sigma-\Psi) H^{-1}$ is the correlation matrix of $\Lambda F$.
Differentiating Eq~(\ref{eq:h_i}) yields
\begin{equation}\label{eq:dh_i}
	dh_i = 2\sum\limits_{t=1}^k \lambda_{it} d\lambda_{it}.
\end{equation}

To obtain the spectral decomposition of this correlation matrix, let $\theta_r$ denote its $r$-th largest eigenvalue, and let $\beta_r$ be the corresponding eigenvector, scaled so that its length is $\theta_r^{1/2}$.
Equivalently, 
$$
\theta_r = \beta_r^\top \beta_r,
$$ 
and 
$$
H^{-1} (\Sigma-\Psi) H^{-1} \beta_r = \theta_r \beta_r, \quad \forall \ 1\le r\le k. 
$$

In alpha factor analysis (\cite{KaiserCaffrey1965}), 
the loading matrix $\Lambda = [\lambda_1,\lambda_2,\ldots,\lambda_k]$ is defined by
$$
\Lambda = H [\beta_1,\beta_2,\ldots,\beta_k]. 
$$
Let $\Gamma = [\gamma_1,\gamma_2,\ldots,\gamma_k]$, where
$\gamma_r = H^{-1} \beta_r$. Then, for each $r=1,\ldots, k$, the
following identities hold:

\begin{equation}\label{eq:alpha-ah2}
\lambda_r = H\beta_r = H^2 \gamma_r,
\end{equation}

\begin{equation}\label{eq:alpha-agh}
\theta_r = \beta_r^\top \beta_r = \gamma_r^\top H^2 \gamma_r = \gamma_r^\top \lambda_r, 
\end{equation}

\begin{equation}\label{eq:alpha-rpslh2g}
\resizebox{.90\textwidth}{!}{$
(\Sigma-\Psi)\gamma_r 
= (\Sigma-\Psi)H^{-1} \beta_r 
= 
H H^{-1} (\Sigma-\Psi)H^{-1} \beta_r
=\theta_r H \beta_r 
= \gamma_r^\top \lambda_r \lambda_r.
$}
\end{equation}
In the remainder of this section, we first use Eq.~(\ref{eq:alpha-ah2}) to express $d\overrightarrow{\Gamma}$ in terms of $d\overrightarrow{\Lambda}$.
We then derive, from Eq.~(\ref{eq:alpha-rpslh2g}), a matrix equation that links $d\overrightarrow{\Sigma}$, $d\overrightarrow{\Gamma}$, and $d\overrightarrow{\Lambda}$.

From Eq.~(\ref{eq:alpha-ah2}), for any $1\le i\le p$ and $1\le r\le k$, we have
$$
\lambda_{ir} = h_i \gamma_{ir}.
$$
Taking differentials yield
$$
d\lambda_{ir} = h_i d\gamma_{ir} + 2\gamma_{ir} \sum_{t=1}^k \lambda_{it} d\lambda_{it}.
$$
Stacking these identities gives the matrix representation
\begin{equation}\label{eq:dGDdA}
	d\overrightarrow{\Gamma} = E\, d\overrightarrow{\Lambda},
\end{equation}
where the entries of $E$ are 
$$
E_{ir,jt} = \frac{\partial \gamma_{ir}}{\partial \lambda_{jt}} 
= 
\delta_{ij}(\delta_{rt}-2\gamma_{ir}\lambda_{it})/h_i,
$$
for all $1\le i, j\le p$ and $1\le r,t\le k$. 

For fixed indices $1\le i \le p$ and $1\le r\le k$, the $i$-th component of Eq.~(\ref{eq:alpha-rpslh2g}) can be written as
\begin{equation}\label{eq:alpha-eq}
\sum_{j\not = i} \sigma_{ij}\gamma_{jr} + h_i \gamma_{ir} = \lambda_{ir} \sum_{j=1}^p \lambda_{jr}\gamma_{jr}.
\end{equation}
Differentiating Eq.~(\ref{eq:alpha-eq}) yields
$$
\resizebox{.98\textwidth}{!}{$
\sum\limits_{j\not = i} \left(\sigma_{ij} d\gamma_{jr} + \gamma_{jr} d\sigma_{ij} \right) 
+ h_i d \gamma_{ir} +  \gamma_{ir} dh_i
=
\lambda_{ir} \sum\limits_{j=1}^p \left(\lambda_{jr} d \gamma_{jr} + \gamma_{jr} d \lambda_{jr}\right)
+ \left(\sum\limits_{j=1}^p \lambda_{jr}\gamma_{jr}\right) d \lambda_{ir}.
$}
$$
Substituting Eq.~(\ref{eq:dh_i}) into the expression above and collecting terms, we obtain
$$
\resizebox{.98\textwidth}{!}{$
\sum\limits_{j\not = i} \gamma_{jr} d\sigma_{ij}
=
\sum\limits_{j\not =i} (\lambda_{ir}\lambda_{jr}-\sigma_{ij}) d\gamma_{jr} +(\lambda_{ir}^2- h_i) d\gamma_{ir} 
+ \left(\sum\limits_{j\not =i}\lambda_{jr} \gamma_{jr}\right) d\lambda_{ir}
- 2\gamma_{ir}\sum\limits_{t\not =r} \lambda_{it} d\lambda_{it} 
+ \lambda_{ir} \sum\limits_{j\not =i} \gamma_{jr}d\lambda_{jr}.
$}
$$

Finally, we collect the differential relations above and write them in matrix form as
\begin{equation}\label{eq:MSNGQA}
F\, d\overrightarrow{\Sigma} = G\, d\overrightarrow{\Gamma} + J\, d\overrightarrow{\Lambda}. 
\end{equation}
The coefficient matrices $F$, $G$, and $J$ are defined as follows.
\begin{itemize}
\item Matrix $F$.
The coefficient on $d\sigma_{xy}$ in the $(i,r)$-th equation of Eq.~(\ref{eq:MSNGQA}) is
$$
F_{ir,xy} = \delta_{ix}(1-\delta_{iy})\gamma_{yr}.
$$
\item Matrix $G$. The coefficient on $d\gamma_{jt}$ is given by
$$
G_{ir,jt} = \left \{
\begin{array}{ll}
\lambda^2_{ir}-h_i,	              		& \mathrm{if} \ j=i \ \mathrm{and} \ t=r,\\
\lambda_{ir}\lambda_{jr}-\sigma_{ij}, 	& \mathrm{if} \ t=r \ \mathrm{but} \ j\not =i,\\
0,				              			& \mathrm{otherwise.}
\end{array}
\right .
$$
\item Matrix $J$. The coefficient on $d\lambda_{jt}$ is
$$
J_{ir,jt} = \left \{
\begin{array}{ll}
\sum\limits_{s\not =i}\lambda_{sr}\gamma_{sr},	& \mathrm{if} \ j=i \ \mathrm{and} \ t=r,\\
-2\gamma_{ir}\lambda_{it},		  				& \mathrm{if} \ j=i \ \mathrm{but} \ t\not =r,\\
\lambda_{ir}\gamma_{jr}, 		  				& \mathrm{if} \ t=r \ \mathrm{but} \ j\not =i,\\
0,				     	          				& \mathrm{otherwise.}
\end{array}
\right .	
$$
\end{itemize}
Combining Eqs.~(\ref{eq:dGDdA}) and (\ref{eq:MSNGQA}) yields
$$
F\, d\overrightarrow{\Sigma} = G\, d\overrightarrow{\Gamma} + J\, d\overrightarrow{\Lambda}
= G\, \left(E\, d\overrightarrow{\Lambda}\right) + J\, d\overrightarrow{\Lambda}
= \left(G\, E + J\right) d \overrightarrow{\Lambda}.
$$
Consequently,
$$
d\overrightarrow{\Lambda} = (G\, E\, + \, J)^{-1}\, F\, d\overrightarrow{\Sigma}.
$$

\section{Image Factor Analysis}\label{sect:ImageFA}
Let $\Delta = \left(\operatorname{Diag} \left(\Sigma^{-1}\right) \right)^{-1}$, where $\operatorname{Diag} \left(\Sigma^{-1}\right)$ denotes the diagonal matrix obtained by retaining only the main diagonal of $\Sigma^{-1}$. 
In image factor analysis (\cite{Joreskog1969}), the covariance matrix $\Sigma$ is modeled as
\begin{equation}\label{eq:se-imageFA}
\Sigma = \Lambda \Lambda^\top + \tau \Delta,
\end{equation}
where $\tau>0$ is an unknown scalar parameter. 
Under this parameterization, one may write $\Psi = \tau \Delta$.
J\"oreskog~(\cite{Joreskog1969}) proposed a maximum-likelihood procedure for estimating this model; in that framework, the covariance matrix of the resulting estimators follows directly from standard maximum-likelihood theory, vulnerable to misspecification of the probability distribution used for $\widehat \Sigma$.

For the principal-factor procedure (\cite{SAS2015}), the loading matrix $\Lambda$ is constructed from
the $k$ largest eigenvalues and their associated eigenvectors of $\Sigma - \tau \Delta$.
The parameter $\tau$ is estimated simultaneously as the slope coefficient from regressing
the diagonal elements of $\Sigma -\Lambda\Lambda^\top$ on the diagonal elements of $\Delta$, namely,
\begin{equation}\label{eq:tau}
\tau = \frac{\langle \operatorname{diag} (\Delta), \operatorname{diag} (\Sigma - \Lambda \Lambda^\top) \rangle}
	{\langle \operatorname{diag} (\Delta), \operatorname{diag} (\Delta) \rangle}
\end{equation}
where $\langle \cdot,\cdot \rangle$ denotes the standard inner product and $\operatorname{diag} (\Delta)$ is the vector collecting the diagonal entries of $\Delta$. 

Let $\sigma^{xy}$ denote the $(x,y)$-th entry of $\Sigma^{-1}$. 
Since $\Sigma \Sigma^{-1}$ is the $p\times p$ identity matrix,
differentiating this identity and applying the product rule yields
$$
d\Sigma^{-1} = - \Sigma^{-1} (d\Sigma) \Sigma^{-1}. 
$$
In particular, for each $i=1,...,p$, the matrix differential identity implies
$$
d\sigma^{ii} = -\sum_{x,y=1}^p \sigma^{ix} \left( d\sigma_{xy} \right) \sigma^{yi}  
= -\sum_{x,y=1}^p \sigma^{ix}\sigma^{iy} d\sigma_{xy}.
$$
Hence,
\begin{equation}\label{eq:d1sigma2}
d\left(\frac{1}{\sigma^{ii}}\right) 
= -\frac{1}{(\sigma^{ii})^2} d\sigma^{ii} 
= \sum_{x,y=1}^p\frac{\sigma^{ix}\sigma^{iy}}{(\sigma^{ii})^2} d\sigma_{xy}.
\end{equation}

Eq.~(\ref{eq:tau}) further implies that 
$$
\tau \langle \operatorname{diag}(\Delta),\operatorname{diag}(\Delta) \rangle
= \langle \operatorname{diag} (\Delta),\operatorname{diag}(\Sigma-\Lambda \Lambda^\top) \rangle, 
$$
and therefore
$$
\tau \sum_{i=1}^p \left( \frac{1}{\sigma^{ii}} \right)^2 = \sum_{i=1}^p \frac{1}{\sigma^{ii}}(\sigma_{ii}-h_i).
$$
Differentiating both sides yields
$$
\resizebox{.98\textwidth}{!}{$
\left[ \sum\limits_{i=1}^p \left( \frac{1}{\sigma^{ii}} \right)^2 \right] d \tau
+ \sum\limits_{i=1}^p \frac{2\tau}{\sigma^{ii}} d \left(\frac{1}{\sigma^{ii}}\right)
=
\sum\limits_{i=1}^p\left[ \frac{1}{\sigma^{ii}}(d \sigma_{ii}-d h_i) + (\sigma_{ii}-h_i) d\left(\frac{1}{\sigma^{ii}}\right)\right].
$}
$$
Substituting Eq.~(\ref{eq:d1sigma2}) into the preceding expression yields
$$
\resizebox{.98\textwidth}{!}{$
	\left[\sum\limits_{i=1}^p \left(\frac{1}{\sigma^{ii}} \right)^2\right] d\tau
= 
\sum\limits_{x,y=1}^p \left[\sum\limits_{i=1}^p \left(\sigma_{ii}-h_i-\frac{2\tau}{\sigma^{ii}} \right)
	\frac{\sigma^{ix}\sigma^{iy}}{(\sigma^{ii})^2}\right] d\sigma_{xy}
	+ \sum\limits_{i=1}^p \frac{1}{\sigma^{ii}} \left(d\sigma_{ii} -2\sum\limits_{t=1}^k \lambda_{it}d\lambda_{it}\right).
$}
$$
These expressions can be written compactly as 
\begin{equation}\label{eq:ifa-se-tbsca}
d\tau = \mu^\top d\overrightarrow{\Sigma} + \eta^\top d\overrightarrow{\Lambda},
\end{equation}
where vectors $\mu\in \mathbb{R}^{p^2}$ and $\eta\in \mathbb{R}^{pk}$. 
In particular, the coefficient multiplying $d\sigma_{xy}$ is 
$$
\mu_{xy}= \frac{
	\frac{\delta_{xy}}{\sigma^{xx}} +\sum\limits_{i=1}^p \left(\sigma_{ii}-h_i-\frac{2\tau}{\sigma^{ii}} \right) 
	\frac{\sigma^{ix}\sigma^{iy}}{(\sigma^{ii})^2} }{
	\sum\limits_{i=1}^p \left(\frac{1}{\sigma^{ii}}\right)^2 }
$$ 
and the coefficient multiplying $d\lambda_{ir}$ is
$$
\eta_{ir} = -\frac{2\lambda_{ir}}{\sigma^{ii}\sum\limits_{z=1}^p \left(\frac{1}{\sigma^{zz}}\right)^2}.
$$

Next, let $\theta_r$ denote the $r$-th largest eigenvalue of $\Sigma - \tau \Delta$, and
let $\lambda_r$ be the corresponding eigenvector scaled to have Euclidean norm $\theta_r^{1/2}$.
Equivalently,  
$$
\lambda_r^\top \lambda_r = \theta_r.
$$ 
Then $\lambda_r$ satisfies the eigenvalue equation
$$
\left(\Sigma-\tau \Delta \right)\lambda_r 
= 
\theta_r \lambda_r = \lambda_r^\top \lambda_r \lambda_r.
$$
For any fixed index $1\le i \le p$, the $i$-th component of this equation is 
$$
\sum_{j=1}^p \sigma_{ij}\lambda_{jr} -\frac{\tau \lambda_{ir}}{\sigma^{ii}} 
= \lambda_{ir} \sum_{j=1}^p \lambda_{jr}^2.
$$
Taking differentials of both sides with respect to the underlying parameters yields
$$
\resizebox{.98\textwidth}{!}{$
\sum\limits_{j=1}^p \left( \sigma_{ij} d \lambda_{jr} + \lambda_{jr} d \sigma_{ij}\right) 
-\frac{\tau}{\sigma^{ii}} d\lambda_{ir} - \frac{\lambda_{ir}}{\sigma^{ii}} d\tau - \tau\lambda_{ir} d \left(\frac{1}{\sigma^{ii}} \right)
= 
2\lambda_{ir} \sum\limits_{j=1}^p \lambda_{jr} d\lambda_{jr} + \left(\sum\limits_{j=1}^p \lambda_{jr}^2\right) d\lambda_{ir}.
$}
$$
Applying Eq.~(\ref{eq:d1sigma2}) and collecting terms involving $d\Sigma, d\Gamma$, and $d\tau$, we obtain
$$
\resizebox{.98\textwidth}{!}{$
	\sum\limits_{j=1}^p \lambda_{jr}d\sigma_{ij}-
	\tau \lambda_{ir}\sum\limits_{x,y=1}^p \frac{\sigma^{ix}\sigma^{iy}}{(\sigma^{ii})^2} d\sigma_{xy}
=
	\sum\limits_{j\not =i}(2 \lambda_{ir} \lambda_{jr}-\sigma_{ij}) d\lambda_{jr} 
	+ \left(2\lambda_{ir}^2-\sigma_{ii}+ \frac{\tau}{\sigma^{ii}} +
	\sum\limits_{j=1}^p \lambda_{jr}^2\right) d\lambda_{ir} + \frac{\lambda_{ir}}{\sigma^{ii}}d\tau.
$}
$$

The above system can be written compactly in matrix form as
\begin{equation}\label{eq:ifa-se-sftga}
	L d\overrightarrow{\Sigma} = P d \overrightarrow{\Lambda} + \pi d\tau
\end{equation}
where $L$, $P$, and $\pi$ are appropriately defined coefficient matrices.
In particular, the coefficient on $d\sigma_{xy}$ in the $(i,r)$-th equation is
$$
L_{ir,xy} = 
\delta_{ix} \lambda_{yr} -\frac{\tau \lambda_{ir}\sigma^{ix}\sigma^{iy}}{(\sigma^{ii})^2}.
$$
The coefficient on $d\lambda_{jt}$ is
$$
P_{ir,jt}= \left \{
\begin{array}{ll}
	2\lambda_{ir}^2-\sigma_{ii}+\frac{\tau}{\sigma^{ii}}+\sum\limits_{z=1}^p \lambda_{zr}^2, 
	& \mathrm{if} \ j=i \ \mathrm{and}\ t=r,\\
	2\lambda_{ir}\lambda_{jr}-\sigma_{ij}, & \mathrm{if} \ t=r \ \mathrm{but}\ j\not = i,\\
	0,                         	& \mathrm{otherwise,}
\end{array}
\right .
$$
and the coefficient on $d\tau$ is
$$
\pi_{ir} = \frac{\lambda_{ir}}{\sigma^{ii}}.
$$
These expressions hold for $1\le i,j,x,y\le p$ and $1\le r,t\le k$. 
Note that $\pi$ collects $\{\pi_{ir} \}$ and thus is a vector of dimension $p k$.

Combining Eqs.~(\ref{eq:ifa-se-tbsca}) and (\ref{eq:ifa-se-sftga}) gives
$$
d\overrightarrow{\Lambda} 
= 
\left(P+ \pi \eta^\top\right)^{-1} \left(L-\pi \mu^\top\right) d\overrightarrow{\Sigma}.
$$ 
Substituting this expression into Eq.~(\ref{eq:ifa-se-tbsca}) yields the differential of $\tau$ in terms of $d\overrightarrow{\Sigma}$:
\begin{equation}\label{eq:tau-sigma}
d\tau 
= 
\left[\mu^\top + \eta^\top \left(P+ \pi \eta^\top\right)^{-1}\left(L-\pi \mu^\top \right)\right] d\overrightarrow{\Sigma}.
\end{equation}
Accordingly, we can also estimate the covariance between $\widehat \Lambda$ and $\widehat \tau$ by
$$
\resizebox{.98\textwidth}{!}{$
\operatorname{cov}\left( \overrightarrow{\widehat{\Lambda}}, \widehat{\tau} \right) 
= 
\left(P+ \pi \eta^\top\right)^{-1} \left(L-\pi \mu^\top\right) 
\operatorname{cov}\left( \overrightarrow{\widehat{\Sigma}} \right) 
\left[\mu^\top + \eta^\top \left(P+ \pi \eta^\top\right)^{-1}\left(L-\pi \mu^\top \right)\right]
$}
$$ 
where $P$, $\pi$, $\eta$, $L$, and $\mu$ are evaluated at $\widehat \Lambda$ and $\widehat \tau$.

\section{Standard Errors of Uniqueness Estimator $\widehat \psi$}\label{sect:uniqueness}
Assume that 
$$
d\overrightarrow{\Lambda} = M\, d\overrightarrow{\Sigma}
$$ 
for some  matrix $M\in \mathbb{R}^{pk\times p^2}$.
To derive the standard errors of the uniqueness estimator $\widehat \psi$, we begin with the identity in Eq.~(\ref{eq:L2}),
$$
\psi_i = \sigma_{ii} - \sum\limits_{t=1}^k \lambda_{it}^2.
$$
Taking differentials on both sides yields
\begin{equation}\label{eq:se-ls-unique}
d\psi_i = d\sigma_{ii} - 2 \sum_{t=1}^k \lambda_{it} d\lambda_{it}, \quad i=1,2,\ldots,p.
\end{equation}
Stacking Eq.~(\ref{eq:se-ls-unique}) across $i$ gives the compact matrix representation
\begin{equation}\label{eq:se-ls-unique-matrix}
d\psi = T d\overrightarrow{\Sigma} + Z d\overrightarrow{\Lambda},
\end{equation}
where the coefficient on $d\sigma_{xy}$ in the $i$-th equation of Eq.~(\ref{eq:se-ls-unique-matrix}) is 
$$
T_{i,xy} = \delta_{ix} \delta_{iy}
$$ 
and the coefficient on $d\lambda_{jt}$ in the same equation is
$$
Z_{i,jt} = -2\delta_{ij} \lambda_{jt},
$$ 
for all $1\le i, j, x, y\le p$ and $1\le t\le k$. 
Substituting $d\overrightarrow{\Lambda} = M\, d\overrightarrow{\Sigma}$ into Eq.~(\ref{eq:se-ls-unique-matrix}) therefore implies
\begin{equation}\label{eq:d_psi_d_sigma}
d \psi = (T+Z M) d\overrightarrow{\Sigma}.
\end{equation}
Consequently,
$$
\operatorname{var}(\widehat \psi) 
= 
(T+Z M) \operatorname{cov}\left(\overrightarrow{\widehat \Sigma}\right) (T+Z M)^\top
$$
and the covariance between $\widehat \Lambda$ and $\widehat \psi$ is estimated by
$$
\operatorname{cov}(\widehat \Lambda, \widehat \psi) 
= 
M \operatorname{cov}\left(\overrightarrow{\widehat \Sigma}\right) (T+Z M)^\top.
$$
When $\Sigma$ is a correlation matrix, $d\sigma_{ii}=0$ for all $i$, so the term $T\, d \overrightarrow{\Sigma}$ drops out and  Eq.~(\ref{eq:d_psi_d_sigma}) simplifies to $d \psi = Z\, M\, d\overrightarrow{\Sigma}$. 

In image factor analysis, the relation
$$
\widehat \psi_i = \widehat \sigma_{ii} - \sum\limits_{t=1}^k \widehat \lambda^2_{it}
$$ 
need not hold. 
We therefore rewrite Eq.~(\ref{eq:tau-sigma}) as 
$$
d\tau = \phi^\top d\overrightarrow{\Sigma}
$$ 
for some vector $\phi \in \mathbb{R}^{p^2}$.
Using the identity $\psi_i = \frac{\tau}{\sigma^{ii}}$ and applying Eq.~(\ref{eq:d1sigma2}), we obtain
\begin{equation}\label{eq:psi-i}
d\psi_i 
=
\frac{d\tau}{\sigma^{ii}} + \tau d\left(\frac{1}{\sigma^{ii}}\right)
=
\frac{1}{\sigma^{ii}}\sum_{x,y=1}^p \phi_{xy} d\sigma_{xy}
	+ \tau \sum_{x,y=1}^p \frac{\sigma^{ix}\sigma^{iy}}{(\sigma^{ii})^2} d\sigma_{xy},
\end{equation}
for all $1\le i\le p$. 
In matrix form, this can be written compactly as 
$$
d\psi = Q\, d\overrightarrow{\Sigma},
$$
where the coefficient multiplying $d\sigma_{xy}$ in the $i$-th equation is
$$
Q_{i,xy} 
= 
\frac{\partial \psi_i}{\partial \sigma_{xy}} 
= 
\frac{\sigma^{ii} \phi_{xy}+\tau\sigma^{ix}\sigma^{iy}}{\left(\sigma^{ii}\right)^2},
\quad \forall \ 1\le i,x,y\le p.
$$

\section{Empirical and Simulation Studies}\label{sect:Example}
This section presents two empirical applications that illustrate how the results derived above can be implemented in practice, together with a simulation study that demonstrates convergence to the corresponding theoretical limits.
We focus in particular on a factor-augmented SVAR application, showing how a factor-embedded model can be used for impulse response analysis, forecast-error variance decomposition, and historical decomposition.

\subsection{Empirical Study: Least-Squares Factorization}
In this subsection, we apply least-squares factor analysis to a sample 
correlation matrix of nine variables measured on $211$ subjects.
The data were reported on page 43 of \cite{Lawley&Maxwell1971}, and the sample is assumed to be drawn from an approximately multivariate normal population.
For simplicity, we set the number of common factors to $k=2$, a choice  
supported by both parallel analysis (\cite{Horn1965}) and the minimum average partial test (\cite{Velicer1976}). 
Estimates of the factor loading matrix $\widehat \Lambda$ and the uniqueness vector $\widehat \psi$ are obtained using the iterative procedure described in \cite{Harman1976}.

To compute the covariances of these estimators, we first apply Eq.~(3.23) in \cite{Girshick1939}
to derive the covariances of the sample correlation coefficients. 
We then use the results developed in Sections~\ref{sect:LSFA} and \ref{sect:uniqueness} to obtain the 
covariance matrix of $\widehat \Lambda$ and the standard errors of $\widehat \psi$. 
To illustrate the impact of factor rotation, we orthogonally rotate $\widehat \Lambda$ using the varimax criterion and compute standard errors for the rotated loadings following the method in \cite{Jennrich&Thayer1973}. 
Table \ref{tb:twofactors} reports the parameter estimates for the least-squares model, with standard errors shown in parentheses.

\begin{table}[H]
\caption{Estimates \& Standard Errors for 2-Factor Least-Squares Factor Analysis}
\label{tb:twofactors}
\centering
\resizebox{.95\hsize}{!}{
\begin{tabular}{c|cc|cc|c}\hline\hline
& \multicolumn{2}{|c|}{Unrotated Factors}& \multicolumn{2}{|c|}{Rotated Factors}&\\
Variable&   I          &   II         & I            & II           &  Uniqueness  \\ \hline
$X_1$   & .6639(.0397) & .3285(.0502) & .6745(.0453) & .3063(.0538) & .4512(.0557) \\
$X_2$   & .6879(.0381) & .2388(.0545) & .6202(.0487) & .3815(.0555) & .4698(.0539) \\
$X_3$   & .4956(.0536) & .2831(.0657) & .5328(.0583) & .2047(.0658) & .6743(.0599) \\
$X_4$   & .8470(.0249) & -.3037(.0372)& .3007(.0381) & .8481(.0295) & .1904(.0412) \\
$X_5$   & .7035(.0392) & -.3179(.0637)& .1990(.0493) & .7459(.0399) & .4040(.0551) \\
$X_6$   & .8037(.0297) & -.3581(.0659)& .2312(.0404) & .8490(.0353) & .2258(.0526) \\
$X_7$   & .6686(.0440) & .3889(.0654) & .7242(.0412) & .2717(.0513) & .4018(.0546) \\
$X_8$   & .4236(.0609) & .2552(.0813) & .4656(.0639) & .1666(.0695) & .7555(.0578) \\
$X_9$   & .7718(.0347) & .4398(.0598) & .8289(.0328) & .3194(.0444) & .2109(.0443) \\ \hline\hline
\end{tabular}
}
\end{table}

\subsection{Simulation Study: Convergence}
To validate the theoretical expression for the covariance $\operatorname{cov}(\widehat \psi)$ derived in Section~\ref{sect:uniqueness}, we conduct a simulation study. 
The sample correlation matrix reported above is treated as the population correlation matrix $\Sigma$.
Using this $\Sigma$, we generate $m$ random correlation matrices from a Wishart distribution with
parameters $(\Sigma, 211)$ for $m=50,100,500,1000,2000$. 
For each simulated correlation matrix, the uniquenesses are estimated by least-squares factoring with $k=2$. 

Based on the resulting set of uniqueness estimates, 
we compute their empirical standard errors.
Table~\ref{tb:simulation} reports these empirical standard errors, with the final column displaying the corresponding theoretical standard errors obtained from the formulas in Section~\ref{sect:uniqueness}. 
Overall, the empirical standard errors converge toward their theoretical counterparts as $m$ increases, although the rate of convergence is relatively slow.

\begin{table}[H]
\caption{Standard Errors for Uniqueness Estimates by Simulations}
\label{tb:simulation}
\centering
\resizebox{.95\hsize}{!}{
\begin{tabular}{c|ccccc|c}\hline\hline
&\multicolumn{5}{|c|}{Empirical S.E. of Uniqueness for Simulated Correlations}&Theoretical\\
Uniqueness& $m=50$   &  $100$   & $500$    & $1000$   & $2000$   &  S.E.     \\ \hline
$\psi_1$  & .0499849 & .0511051 & .0549498 & .0560314 & .0556181 & .0556623  \\
$\psi_2$  & .0513352 & .0520701 & .0545346 & .0555878 & .0549598 & .0538818  \\
$\psi_3$  & .0666649 & .0660852 & .0599577 & .0603062 & .0598762 & .0598473  \\
$\psi_4$  & .0405467 & .0433619 & .0432832 & .0428048 & .0415737 & .0411420  \\
$\psi_5$  & .0481342 & .0538724 & .0524052 & .0530417 & .0544235 & .0550669  \\
$\psi_6$  & .0501591 & .0543790 & .0541692 & .0550827 & .0546638 & .0526150  \\
$\psi_7$  & .0524787 & .0565299 & .0547573 & .0553479 & .0546232 & .0546297  \\
$\psi_8$  & .0507511 & .0562004 & .0576808 & .0582864 & .0580424 & .0578393  \\
$\psi_9$  & .0394781 & .0409714 & .0430909 & .0448280 & .0438135 & .0443326  \\ \hline\hline
\end{tabular}
}
\end{table}

\subsection{Empirical Study: Factor-Augmented SVAR}
We consider the SVAR model specified in Eq.~(\ref{eq:SVAR}) using three endogenous time series: potential GDP growth (gdp), the output gap (gap), and the unemployment rate (unrate).
The data are obtained from the Federal Reserve Economic Database (FRED, \url{https://fred.stlouisfed.org/}) and cover the period 2000--2025.

The output gap is typically viewed as predominantly demand-driven, particularly when potential (trend) output is interpreted as reflecting supply-side conditions.
Accordingly, for $Y_t = (gdp_t, \;gap_t, \;unrate_t)^\top$, we impose the following zero restrictions on the loading matrix:
\begin{equation}\label{eq:restrictions}
\Lambda = \left[
\begin{array}{cc}
NA & 0 \\
0  & NA \\
NA & NA
\end{array}\right]
\end{equation}
where $F_t$ collects the supply and demand shocks, in that order,
and all remaining shocks are absorbed into $\epsilon_t$.

Given the pattern matrix specified in Eq.~(\ref{eq:restrictions}), the model in Eq.~(\ref{eq:model}) is still not exactly identified by the zero restrictions alone, since the restricted $\Lambda$ contains four unknowns while $\Sigma$ has only three distinct off-diagonal elements.
We therefore choose $\Lambda$ so that $\Lambda \Lambda^\top$ reproduces the off-diagonal elements of $\Sigma$ and, subject to this constraint, fits the diagonal elements of $\Sigma$ as closely as possibly---thereby maximizing the explanatory power of the supply and demand shocks.
This leads to the constrained minimization problem:
\begin{equation}\label{eq:restricted_min}
\min\limits_{_{0\le \lambda_{11} \le \sqrt{\sigma_{11}},\; 0\le \lambda_{22} \le \sqrt{\sigma_{22}}}}
\left[ 
\left(\sigma_{11}-\lambda_{11}^2\right) 
+ \left(\sigma_{22}-\lambda_{22}^2\right) 
+ \left(\sigma_{33} - \lambda_{31}^2-\lambda_{32}^2\right)
\right] 
\end{equation}
subject to 
$$
\lambda_{11} \lambda_{31} = \sigma_{31}, \quad
\lambda_{22} \lambda_{32} = \sigma_{32}, \quad
0\le \lambda_{31}^2 + \lambda_{32}^2 \le \sigma_{33}.
$$

Although the first-order condition associated with Eq.~(\ref{eq:restricted_min}), or the Lagrange multiplier of the constrained minimization problem, can in principle serve as an identifying restriction $\Phi(\Lambda, \Sigma)=0$,
implicit differentiation may perform poorly in this setting because the quadratic minimization is frequently resolved at the boundary.
In practice, however, standard numerical algorithms can readily compute $\widehat \Lambda$ for the specification in Eq.~(\ref{eq:restricted_min}). 

To estimate the common factors $F_t$ in Eq.~(\ref{eq:SVAR}), we first construct the estimated two-component reduced-form residuals,
$$
\widehat u_t 
= 
Y_t - \widehat c - \sum_{s=1}^{\zeta} \widehat A_s Y_{t-s},
$$
and then recover $\widehat F_t$ from 
$$
\widehat u_t = \widehat \Lambda \widehat F_t + \widehat \epsilon_t 
$$ 
using factor-score methods (e.g., \cite{Harman1976}).
The estimated idiosyncratic shocks are subsequently obtained as
$$
\widehat \epsilon_t = \widehat u_t - \widehat \Lambda \widehat F_t.
$$

\begin{figure}[H]
\centering
\includegraphics[height=8cm, width=13cm]{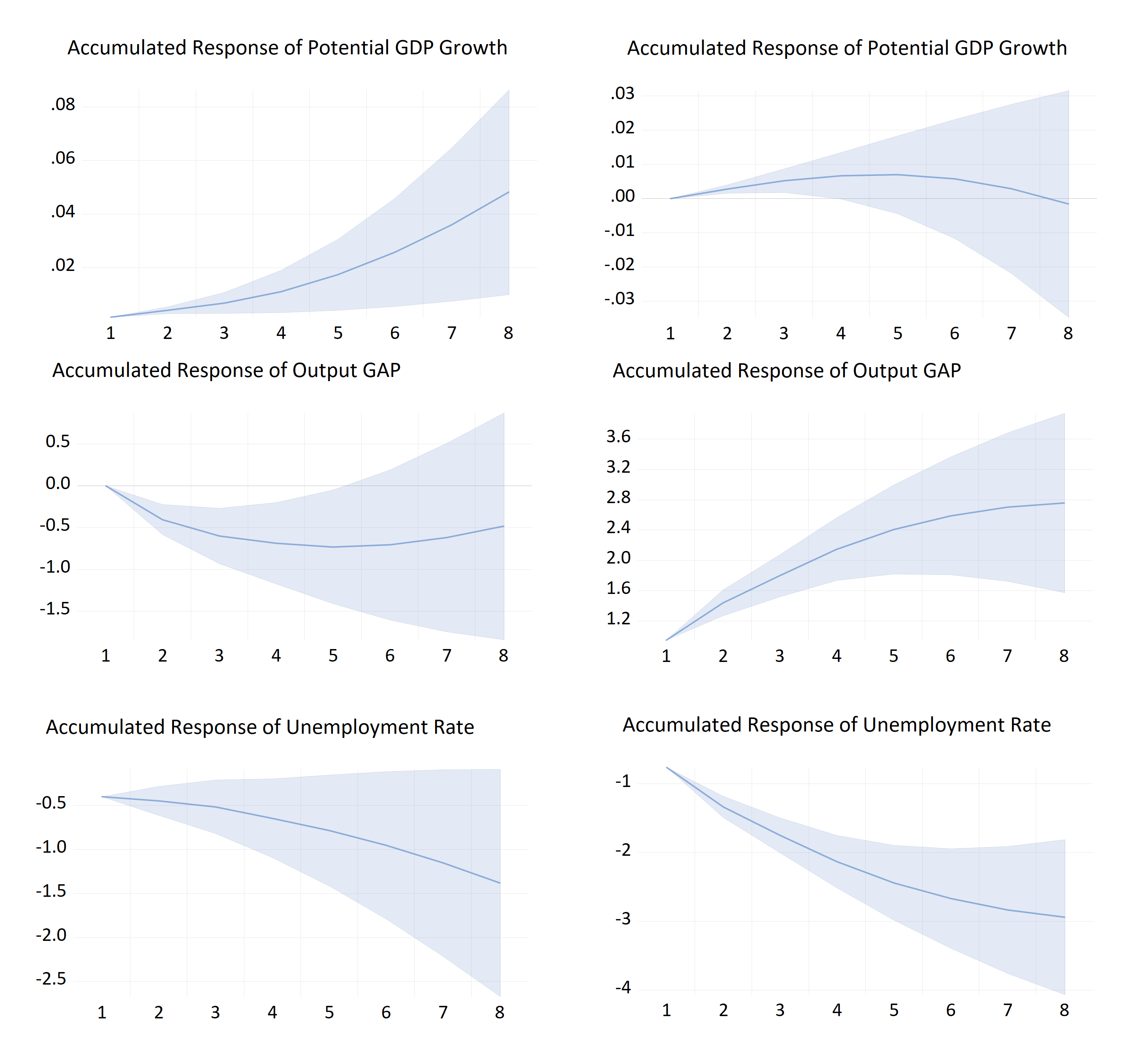}
\caption{Responses to Supply (left) and Demand (right) Shocks when $\operatorname{cov}(\widehat \Lambda)$ is Ignored}
\label{Fig:IRF_wo_cov}
\end{figure}

For the impulse-response analysis, a one-unit supply (demand) shock is introduced into the SVAR system through the first (second) column of $\widehat \Lambda$. 
Figure~\ref{Fig:IRF_wo_cov} reports the responses of $Y_t$ to unit common shocks, treating $\operatorname{cov}(\widehat \Lambda)$ as negligible (i.e., ignored). 
Although the responses exhibit the expected signs, their dynamic paths differ markedly across the two shocks, with notable differences in both curvature and magnitude. 
Moreover, both common shocks---as measured at horizon $1$---display no uncertainty.
The $95\%$ confidence intervals beyond horizon 1 reflect only the sampling variability in the estimated coefficient matrices.

\begin{figure}[H]
\centering
\includegraphics[height=8cm, width=13cm]{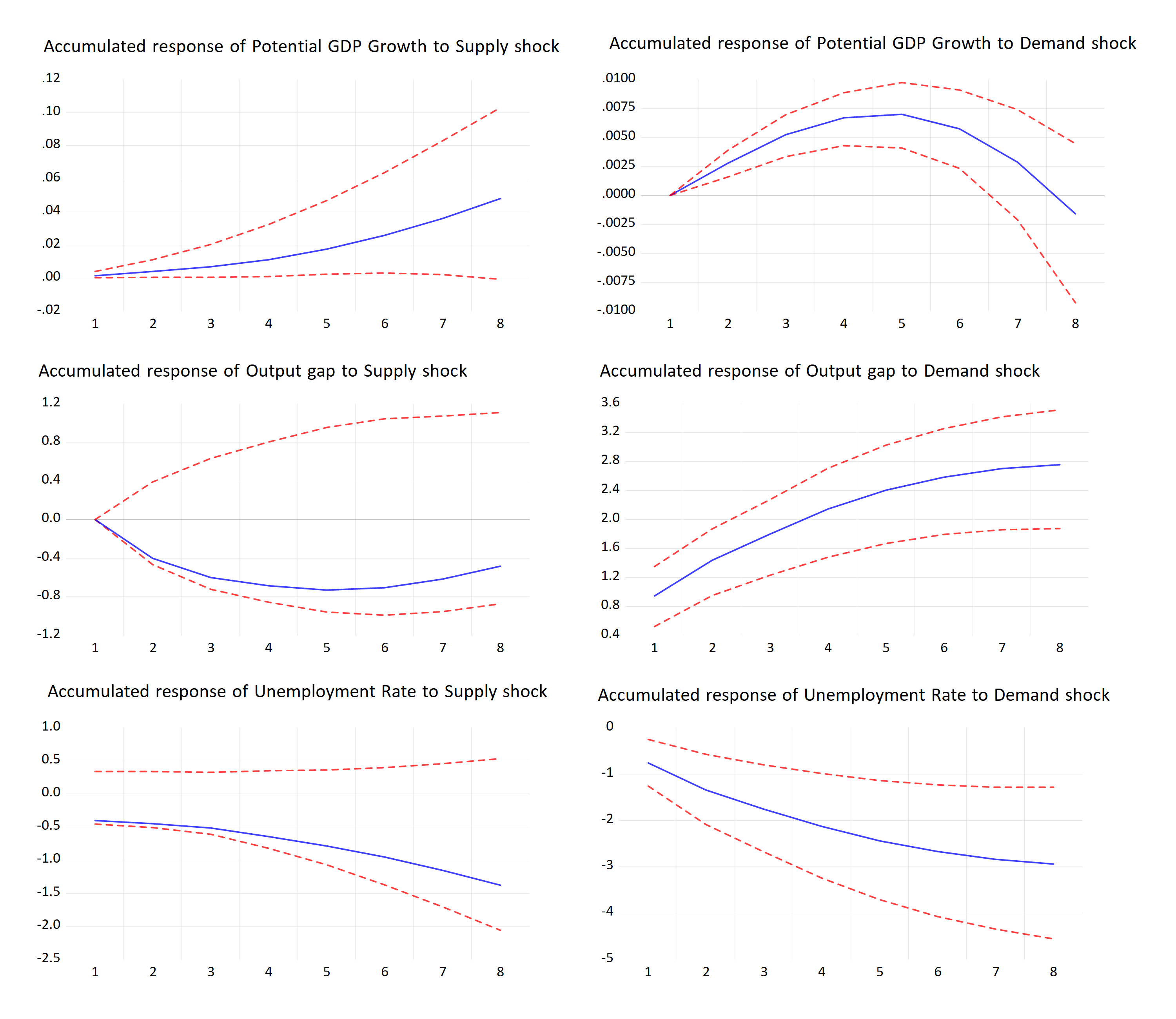}
\caption{Responses to a Unit Supply or Demand Shock when $\operatorname{cov}(\widehat \Lambda)$ is Incorporated}\label{Fig:IRF_w_cov}
\end{figure}

To incorporate $\operatorname{cov}(\widehat \Lambda)$ into the impulse-response analysis, we employ a residual bootstrap that generates $5,000$ bootstrap replications of the variance-covariance matrix $\widehat \Sigma$.
For each replication, the corresponding $\widehat \Lambda$ is computed via Eq.~(\ref{eq:restricted_min}), yielding $5,000$ bootstrap impulse-response functions.
Figure~\ref{Fig:IRF_w_cov} reports the resulting 95\% confidence intervals for the responses of $Y_t$ to unit common shocks once the uncertainty in $\widehat \Lambda$ is taken into account. 
Shock uncertainty is summarized by the confidence intervals at horizon $1$, where
exactly two responses exhibit no uncertainty as a direct consequence of the identifying restrictions imposed in Eq.~(\ref{eq:restrictions}).

A comparison of the impulse responses computed with and without $\operatorname{cov}(\widehat \Lambda)$ reveals that $\operatorname{cov}(\widehat \Lambda)$ is far from negligible.
As shown in Figure~\ref{Fig:IRF_w_cov}, incorporating $\operatorname{cov}(\widehat \Lambda)$ materially widens the confidence bands---with one exception---and the nonlinearity inherent in the factor extraction induces pronounced asymmetry.
The confidence intervals in Figure~\ref{Fig:IRF_w_cov} therefore reflect the joint uncertainty in both $\widehat A_s$ and $\widehat \Lambda$.

\subsection{Historical and Forecast-Error Variance Decomposition}
Historical decomposition expresses $Y_t$ as the sum of mutually independent components driven by the two common shocks in $F_t$ and the three idiosyncratic shocks in $\epsilon_t$.
Rewriting Eq.~(\ref{eq:SVAR}) in its moving-averaging (Wold) representation replaces lagged values of $Y_t$ with lagged shocks, yielding
\begin{equation}\label{eq:MA_Representation}
Y_t 
= 
\tilde{c} 
+
\sum\limits_{l=0}^\infty \omega_l F_{t-l}
+ 
\sum\limits_{l=0}^\infty \Omega_l \epsilon_{t-l}
=
\tilde{c}
+
\sum\limits_{l=0}^\infty \sum\limits_{i=1}^k \omega_l F_{t-l}^{(i)}
+ 
\sum\limits_{l=0}^\infty \sum\limits_{j=1}^p \Omega_l \epsilon_{t-l}^{(j)}.
\end{equation}
Here, $F_{t-l}^{(i)}$ denotes the vector that retains only the $i$-th element of $F_{t-l}$, with all other entries set to zero;
analogously, $\epsilon_{t-l}^{(j)}$ is the vector containing only the $j$-th element of $\epsilon_{t-l}$.
The coefficient matrices $\omega_l$ and $\Omega_l$ are polynomials in $A_s$ and $\Lambda$;
in particular, $\omega_l$ corresponds to the impulse-response coefficient matrix of $Y_t$ with respect to the common shocks $F_{t-l}$.

Accordingly, the contribution of the $i$-th common shock to $Y_t$ is 
\begin{equation}\label{eq:histdecomp_common}
\sum\limits_{l=0}^\infty \omega_l F_{t-l}^{(i)},
\end{equation}
and the contribution of the $j$-th idiosyncratic shock is 
\begin{equation}\label{eq:histdecomp_idiosyn}
\sum\limits_{l=0}^\infty \Omega_l \epsilon_{t-l}^{(j)}.
\end{equation}
In estimation, the quantities $F_t$, $\epsilon_t$, $\omega_l$ and $\Omega_l$ are replaced by their corresponding sample counterparts. 
Both $\widehat F_t$ and $\widehat \epsilon_t$ inherit estimation uncertainty from $\widehat \Sigma$,
and their values also depend on the particular factor-score method employed.

For the forecast-error variance decomposition of $Y_{t+h}$ (forecast at time $t$), Eq.~(\ref{eq:MA_Representation}) implies that the conditional variance-covariance matrix of $Y_{t+h}$ is
$$
\operatorname{cov}\left(Y_{t+h}|Y_t\right)
= 
\sum\limits_{l=0}^h \omega_l \omega_l^\top
+ 
\sum\limits_{l=0}^h \Omega_l \Psi \Omega_l^\top.
$$ 
Within this decomposition, Eq.~(\ref{eq:histdecomp_common}) implies that the contribution of the $i$-th common shock is 
\begin{equation}\label{eq:FEVD_common}
\sum\limits_{l=0}^h \omega_l Z^{(i)} \omega_l^\top,
\end{equation}
where $Z^{(i)}$ is the $k\times k$ selection matrix with a single one in the $(i,i)$-th position and zeros elsewhere.
Analogously, Eq.~(\ref{eq:histdecomp_idiosyn}) implies that the contribution of the $j$-th idiosyncratic shock to $\operatorname{cov}\left(Y_{t+h}|Y_t\right)$ is
\begin{equation}\label{eq:FEVD_idio}
\sum\limits_{l=0}^h \Omega_l \Psi^{(j)} \Omega_l^\top,
\end{equation}
where $\Psi^{(j)}$ is the $p\times p$ matrix that retains only the $(j,j)$-th entry of $\Psi$, with all other entries set to zero.

When $\widehat \Omega_l$ and $\widehat \Psi$ are substituted into Eq.~(\ref{eq:FEVD_idio}), the two estimators are generally correlated, since both $\operatorname{cov}(\widehat \Omega_l)$ and $\operatorname{cov}(\widehat \Psi )$ ultimately derive from $\operatorname{cov}(\widehat \Sigma )$.
For this reason, a resampling-based procedure is typically preferable to a purely analytical approach for computing the standard errors associated with Eq.~(\ref{eq:FEVD_idio}), whereas either method remain feasible for Eq.~(\ref{eq:FEVD_common}).

\section{Conclusion and Discussions}\label{sect:diss}

Recent years have seen renewed interest in formal statistical inference for factor-analytic models. 
This resurgence reflects several practical developments, including the increasing prevalence of non-normal data, heightened concerns about robustness, more complex covariance structures,
and the expansion of both cross-sectional and temporal dimensions~(\cite{Bai2003,Fan&Wang&Zhong&Zhu2021}).
Yet, explicit covariance formulas for non-ML extraction methods remain comparatively scarce, despite the continued prevalence of these methods in empirical research and their implementation in standard statistical software (\cite{Harman1976,SAS2015}). 
This paper helps close this gap by deriving analytically tractable covariance expressions for several non-ML factor solutions. 
By embedding these results within a unified implicit-differentiation framework, the proposed approach accommodates non-normal or non-i.i.d. data and integrates naturally with existing delta-method results for factor rotation (\cite{Archer&Jennrich1973,Jennrich1973,Yung&Hayashi2001}) and other factor solutions. 
The framework enables transparent computation of standard errors without recourse to resampling
and complements recent advances in robust and high-dimensional factor modeling by providing a rigorous inferential foundation for classical extraction methods.

The availability of closed-form covariance also provides a useful diagnostic tool for assessing and validating model specification. 
Unusually large or otherwise anomalous estimated standard errors in $\widehat \Lambda$ and $\widehat \psi$ may indicate model misspecification or systematic estimation errors rather than genuine sampling variability. 
Such distortions can arise from over-factoring, when $k$ is chosen too large, or under-factoring, when $k$ is chosen too small. 
They may also reflect a weak factor structure, an inappropriate extraction procedure, or insufficient separation between  common and idiosyncratic components. 
In particular, if a large fraction of the estimated coefficients in $\widehat \Lambda$ is statistically insignificant, the pattern may be indicative of over-factoring.
Conversely, an overwhelmingly large share of strong and statistically significant coefficients may suggest under-factoring.

Beyond inference on structural shocks, the proposed covariance formulas also facilitate uncertainty assessment in empirical workflows where factor extraction serves as an intermediate step rather than the final objective. 
In applied macroeconomics and finance, factor-analytic structures arise routinely in dynamic factor models, large Bayesian or classical VARs, factor-augmented regressions, and error-component representations of high-dimensional systems. 
Such applications typically involve multiple stages, yet closed-form expressions for the uncertainty of downstream quantities---such as forecasts or policy counterfactuals---are rarely available.
By applying the chain rule, our results provide analytic formulas that propagate uncertainty in $\widehat \Sigma$ through the entire estimation pipeline and quantify the resulting uncertainty in forecasts, impulse responses, or policy conclusions.
In the context of factor rotation, for example, several criteria, such as varimax, quartimax, and promax, are available; an effective rotation criterion should attenuate, rather than amplify, uncertainty when transforming unrotated factors into rotated factors. 

For the SVAR model in Eq.~(\ref{eq:SVAR}), target or Procrustean rotation is unnecessary because the structural shocks $F_t$ are already aligned with the endogenous variables $Y_t$ through the loading matrix.
This alignment naturally induces identifying restrictions---such as zero or sign restrictions---on elements of $\Lambda$. 
Unlike the restrictions commonly used in the SVAR literature, the dimension of the structural factors in this setting is more tractable. 
Nevertheless, the covariance formulas developed in this paper remain applicable because they treat the estimated $\Lambda$ as an input. 
They also remain valid under heteroskedasticity, serial dependence, and non-Gaussian sampling, features that are pervasive in macroeconomic and financial data. 
The framework therefore integrates naturally with HC and HAC estimators while avoiding the strong likelihood-based assumptions that are often difficult to justify in practice.
For example, when the data are multivariate binary or categorical---settings frequently encountered in applied factor analysis---computing $\operatorname{cov}(\widehat \Lambda)$ poses no substantive practical difficulty.

For expositional simplicity, we do not explicitly impose the symmetry conditions $\sigma_{ij} = \sigma_{ji}$ and $\widehat \sigma_{ij} =\widehat \sigma_{ji}$, nor do we spell out the associated relationships among the remaining elements of $\Sigma$ and $\widehat \Sigma$.
In practice, however, exploiting symmetry is important for reducing memory requirements and improving computational efficiency, particularly in large-scale applications.
Moreover, reducing the dimension of $d \Sigma$ from $p^2$ to $p(p+1)/2$ can help mitigate near-singularity issues that may arise when inverting the matrices appearing in these formulas.

Finally, the proposed covariance formulas extend naturally to econometric settings involving endogeneity.
In such applications, the relevant $\widehat \Sigma$ and $\mathrm{cov}(\widehat \Sigma)$ are constructed from second-stage residuals, for example, those obtained from instrumental-variable, generalized method of moments, or control-function estimation within a system of equations. 
The key requirement is the availability of consistent estimators of $\widehat \Sigma$ and $\operatorname{cov}(\widehat \Sigma)$, such as HAC or HC estimators, as implied by the chosen endogeneity-robust procedure.
The resulting $\mathrm{cov}(\widehat \Lambda)$ can then be used to conduct inference on model specification and other functionals of the estimated factor structure.
For certain downstream analyses, such as impulse-response function, the common or idiosyncratic factors should be specified so as to explicitly account for endogeneity.



\appendix
\section{Abbreviations Used in This Manuscript}

\resizebox{.95\hsize}{!}{
\begin{tabular}{ll}
$Y$&Vector of $p$ correlated random variables\\
$\Sigma$&Unobservable population variance-covariance matrix of $Y$\\
$\widehat \Sigma$&Sample variance-covariance matrix of the random vector $Y$\\
$\Lambda$&Loading matrix in the variance-covariance structure model $\Sigma = \Lambda \Lambda^\top +\Psi$ \\
$\Psi$&Uniqueness (a diagonal matrix) of idiosyncratic residuals\\
$\widehat \Lambda$, $\widehat \Psi$&Estimates of the loading matrix and uniqueness, respectively\\
$(p,k)$&Dimensions of $\Lambda$; $p$ is also the dimension for $Y$ and $p\times p$ for $\Sigma$ and $\Psi$\\
$\delta_{ij}$&The Kronecker delta, that equals $1$ when $i=j$ and $0$ otherwise\\
$\overrightarrow{\Lambda}$&The vectorization of the matrix $\Lambda$, ordered by columns. Similar for $\overrightarrow{\Sigma}$\\
$\lambda_{ir}$&The $(i,r)$-th entry of $\Lambda$ with $1\le i\le p$ and $1\le r\le k$\\
$\psi$&The main diagonal of $\Psi$---i.e., the vector of diagonal elements of $\Psi$\\
$\theta_r$&The $r$-th largest eigenvalue of a symmetric matrix\\
$\sigma_{xy}$&The $(x,y)$-th entry of $\Sigma$\\
$\sigma^{xy}$&The $(x,y)$-th entry of $\Sigma^{-1}$---the inverse of $\Sigma$\\
$\operatorname{diag}(M)$&The vector of diagonal elements of the square matrix $M$\\
$\operatorname{Diag}(M)$&The diagonal matrix with diagonal elements in $M$\\
cov($\widehat \Sigma$)&Variance-covariance matrix of $\widehat \Sigma$\\
SVAR&Structural Vector AutoRegressive model\\
$\mathbb{R}$&The set of real numbers. Thus, $\mathbb{R}^{p\times k}$ is the set of $p\times k$ matrices\\
$\omega_l$, $\Omega_l$&Coefficient matrices for the MA or Wold representation of time series\\
$\Phi(\Lambda,\Sigma)$&$\Phi(\Lambda,\Sigma)=0$ is the identifying condition for a factor extraction method\\
$u_t$&$u_t=\Lambda F_t +\epsilon_t$ is the two-component residual in a factor-augmented SVAR model\\
ML&maximum-likelihood\\
non-ML&Non-maximum-likelihood\\
HAC&Heteroskedasticity- and Autocorrelation-Consistent\\
HC&Heteroskedasticity-Consistent
\end{tabular}
}

\vskip 1cm
\noindent \textbf{Declaration of competing interests}: None

\noindent \textbf{Declaration of generative AI use}: None

\noindent \textbf{Funding sources}: None
\end{document}